\documentclass[11pt]{article}
\usepackage{amsfonts}
\usepackage{amsmath,amssymb}

\newtheorem{theorem}{Theorem}
\newtheorem{proposition}[theorem]{Proposition}

\begin{document}
\large

\begin{center}
\Large\bf QUANTUM GROUPS AND NON-COMMUTATIVE COMPLEX ANALYSIS
\end{center}

\bigskip

\centerline{\tt S. Sinel'shchikov and L. Vaksman}

\bigskip

The problem of uniform approximation by analytic polynomials on a compact
$K\subset\mathbb{C}$ made an essential impact to the function theory and
the theory of commutative Banach algebras. This problem was solved by S.
Mergelyan. Much later a theory of uniform algebras was developed and an
abstract proof of Mergelyan's theorem was obtained.

The work by W. Arveson \cite{Arv1} starts an investigation of
non-commutative analogs for uniform algebras. In particular, a notion of the
Shilov boundary for a subalgebra of a $C^*$-algebra have been introduced
therein. So, the initial results of non-commutative complex analysis were
obtained. We assume the basic concepts of that work known to the reader.

In mid'90-s an investigation of quantum analogs for bounded symmetric
domains has been started within the framework of the quantum group theory
\cite{Drinf1}. The simplest of those is a unit ball in $\mathbb{C}^n$. Our
goal is to explain that the quantum sphere is the Shilov boundary for this
quantum domain. The subsequent results of the authors on non-commutative
function theory and quantum groups are available at www.arxiv.org .
Specifically, we obtained some results on weighted Bergman spaces, the
Berezin transform, and the Cauchy-Szeg\"o kernels for quantum bounded
symmetric domains introduced in \cite{SV1}.

In what follows the complex numbers are assumed as a ground field and all
the algebras are assumed to be unital. In what follows $q\in(0,1)$.

To introduce a quantum unit ball, consider a $*$-algebra
$\operatorname{Pol}(\mathbb{C}^n)_q$ given by the generators
$z_1,z_2,\ldots,z_n$ and the defining relations $z_jz_k=qz_kz_j$ for $j<k$,
\begin{eqnarray*}
z_j^*z_k=qz_kz_j^*,\quad j\ne k,\qquad
z_j^*z_j=q^2z_jz_j^*+(1-q^2)\left(1-\sum_{k>j}z_kz_k^*\right).
\end{eqnarray*}
This $*$-algebra has been introduced by W. Pusz and S. Woronowicz
\cite{PWor} where one can find a description (up to unitary equivalence) of
its irreducible $*$-representations $T$. One can demonstrate that
$0<\|f\|\overset{\mathrm{def}}{=}\sup\limits_T\|T(f)\|<\infty$ for all
non-zero $f\in\operatorname{Pol}(\mathbb{C}^n)_q$ and that its
$C^*$-enveloping algebra $C(\mathbb{B})_q$ is a q-analogue for the
$C^*$-algebra of continuous functions in the closed unit ball in
$\mathbb{C}^n$. A plausible description of this $C^*$-algebra has been
obtained by D.~Proskurin and Yu.~Samoilenko \cite{PrSam1}.

To introduce a quantum unit sphere, consider a closed two-sided ideal $J$ of
the $C^*$-algebra $C(\mathbb{B})_q$ generated by
$1-\sum\limits_{j=1}^nz_jz_j^*$. Obviously, the $C^*$-algebra
$C(\partial\mathbb{B})_q\overset{\mathrm{def}}{=}C(\mathbb{B})_q/J$ is a
q-analogue for the algebra we need. Thus the canonical onto morphism
$$j_q:C(\mathbb{B})_q\to C(\partial\mathbb{B})_q$$
is a $q$-analogue for the restriction operator of a continuous function onto
the boundary of the ball.

The closed subalgebra $A(\mathbb{B})_q\subset C(\partial\mathbb{B})_q$
generated by $z_1,z_2,\ldots,z_n$ is a $q$-analogue for the algebra of
continuous functions in the closed ball which are holomorphic in its
interior.

Let $j_{A(\mathbb{B})_q}$ be the restriction of the homomorphism $j_q$ onto
the subalgebra $A(\mathbb{B})_q$.

\begin{theorem}
The homomorphism $j_{A(\mathbb{B})_q}$ is completely isometric.
\end{theorem}

This result is a q-analogue of the well known maximum principle for
holomorphic functions. By an Arveson's definition \cite{Arv1} this means
that $J$ is a boundary ideal for the subalgebra $A(\mathbb{B})_q$. A proof
of this theorem elaborates the methods of quantum group theory and theory of
unitary dilations \cite{SeNaFo}.

One can use the deep result of M. Hamana \cite{Ham1} (on existence of the
Shilov boundary) to prove the following simple

\begin{proposition}
$J$ is the largest boundary ideal for $A(\mathbb{B})_q$.
\end{proposition}

Thus the quantum sphere is the Shilov boundary of the quantum ball. A more
detailed exposition of this talk is available in \cite{Vaks}.

The second named author would like to express his gratitude to
M.~Liv\v{s}itz and V.~Drinfeld who taught him the theory of non-selfadjoint
linear operators and quantum groups, respectively.

\end{document}